\documentclass[12pt, a4paper]{article}
\usepackage[latin1]{inputenc}
\usepackage[english]{babel}
\usepackage{indentfirst}
\usepackage{amstext}
\usepackage{setspace}
\usepackage{amsfonts}
\usepackage{textcomp}
\usepackage{amssymb}
\usepackage{amscd}
\usepackage{epsf}
\usepackage{graphicx}
\usepackage{epsfig}

\newcommand {\demo}{\hskip -0.6cm{\bf Proof:  }}
\newcommand {\fim}{\hfill{$\square$}\vskip 1pc}

\newcommand {\N}{\mathbb{N}}
\newcommand {\Z}{\mathbb{Z}}

\newcommand {\C}{\mathbb{C}}

\newtheorem{theorem}{Theorem}[section]
\newtheorem{lema}[theorem]{Lemma}

\newtheorem{definition}[theorem]{Definition}
\newtheorem{proposition}[theorem]{Proposition}
\newtheorem{example}[theorem]{Example}
\newtheorem{obs}[theorem]{Remark}

\begin{document}

\title{C*-algebras Associated to Stationary Ordered Bratteli Diagrams}
\maketitle
\begin{center}
{\large Daniel Gonçalves and Danilo Royer}\\
\end{center}  
\vspace{8mm}

\abstract 
In this paper, we introduce a C*-algebra associated to any substitution (via its Bratteli diagram model). We show that this C*-algebra contains the partial crossed product C*-algebra of the corresponding Bratteli-Vershik system and show that these algebras are invariant under equivalence of the Bratteli diagrams. We also show that the isomorphism class of the algebras, together with a distinguished set of generators, is a complete invariant for equivalence of Bratteli diagrams.

\onehalfspace

\section{Introduction}



Bratteli-Vershik systems provide an interesting combinatorial model for many dynamical systems, as it is the case of substitutional dynamical systems (see \cite{DHS}) and Cantor minimal systems (see \cite{HPS}) and are also related to cellular automata (see \cite{PY}). From the non-commutative geometry point of view, the object usually associated to a Bratteli-Vershik system is the crossed product C*-algebra. Recently, Fujino introduced another C*-algebra associated to a Bratteli-Vershik system arising from a primitive substitution, see \cite{fujino}. But, as far as we know, the literature is missing a non commutative counterpart that can be defined for any substitution, or any ordered (not necessarily well ordered) stationary Bratteli diagram. It is with this goal that we write this paper. 

Building from the the ideas in \cite{fujino} we define a C*-algebra that can be constructed from any stationary ordered Bratteli diagram (and hence from any substitution) and that has a quotient that is isomorphic to Fujino's algebra (in the case of diagrams arising from primitive, proper substitutions). Of course one of the hurdles we have to overcome when dealing with non necessarily well ordered Bratteli diagrams is the fact that there may be more than one maximal and one minimal infinite path in the diagram, and so it is not clear how to extend the Vershik map to the whole path space. We overcome this problem considering the Vershik map as a partially defined map, and define the C*-algebra accordingly (see section 2). 

One of the key features of Fujino's algebra is that it contains the crossed product C*-algebra arising from the Bratteli-Vershik system. We show, in section 3, that our algebra contains the partial crossed product arising from the Bratteli-Vershik system and take the opportunity to show that this partial crossed product is an AF algebra. Finally, in section 4, we show that the new algebras introduced are invariant under equivalence of the Bratteli diagrams and show that the isomorphism class of the algebras, together with a distinguished set of generators, is a complete invariant of equivalence, that is, we show that if two diagrams are equivalent their algebras are isomorphic (via an isomorphism that preserves the generators), and if two diagrams have algebras that are isomorphic, via an isomorphism that preserves the generators of these algebras, then the diagrams are equivalent.

Before we proceed, we recall the relevant definitions about Bratteli diagrams (following \cite{DHS}):

A Bratteli diagram, (V,E), is an infinite directed graph. The vertex set V is the union of a sequence of finite, non-empty, pairwise, disjoint sets, $V_n$, $n\geq 0$. The set $V_0$ is assumed for convenience to consist of a
single vertex, $v_0$. Similarly, the edge set is the union of a sequence of finite, non-empty, pairwise disjoint sets, $E_n$, $n\geq 0$. An edge $e$ in $E_n$ has initial vertex $i(e)$ in $V_n$ and terminal vertex $t(e)$ in $V_{n+1}$. The graph is always assumed to have no sources other than $v_0$ and no sinks, that is, $i^{-1}(v)$ and
$t^{-1}(v)$ are non-empty for any $v$ in $V$ (other than $t^{-1}(v_0)$). 

A Stationary Bratteli diagram is a Bratteli diagram, $(V,E)$, such that each $V_n$, $n\geq 1$, consists of $N$ points, say $$V_n=\left\{x(1,n), x(2,n),\ldots, x(N,n)\right\},$$ and each $E_n$, $n\geq 1$ consists of $M$ edges. Following \cite{fujino}, $j$ is said to be the label of the vertex $x(j,n)$ and is denoted by $L(x(j,n))$. Furthermore, the edges at level one completely determine the diagram, that is, for each edge $e'$ in $E_1$ with $i(e')=x(i,1)$ and $t(e')=x(j,2)$ there exists an edge $e$ in $E_n$, $n\geq 2$, with $i(e)=x(i,n)$ and $t(e)=x(j,n+1)$ and all edges in $E_n$ arise in this way.

Throughout this paper it is assumed that there is a partial order on $E_n$, $n\geq 1$, such that two edges are comparable if an only if they have the same terminal vertex. A stationary Bratteli diagram is said ordered, if it has `` the same" order on each edge level $E_n$, that is, all $E_n$, $n\geq 1$ are order isomorphic. The set of maximal edges in $E_n$ is denoted by $E_n^{\max}$ and analogously $E_n^{\min}$ denotes the set of minimal edges. Finally, the successor map (Vershik map) on $E_n$, $n\geq 1$, is denoted by $\lambda$.

\section{The C*-algebra associated to an ordered stationary Bratteli diagram}
\label{section2}

In \cite{fujino}, Fujino defined a C*-algebra associated to a primitive substitution. This C*-algebra was built as an universal C*-algebra generated by partial isometries and an unitary (implementing the Vershik map) subject to relations. For primitive, proper substitutions, the associated Bratteli diagrams have only one maximal and one minimal infinite path, and so the Vershik map may be extended to the whole infinite path space. For a general substitution this is not the case and it is not clear how to extend the Vershik map. So the natural idea to overcome this problem is to consider the Vershik map as a partially defined map, which leads to a partial isometry in the universal C*-algebra, instead of an unitary. This is the key point of our definition of the C*-algebra associated to an ordered stationary Bratteli diagram:

\begin{definition}\label{relations} Let $(V,E,\leq)$ be a stationary ordered Bratteli diagram. Define $B_E$ to be the universal unital C*-algebra generated by partial isometries $\{s_e \}_{e\in E_1}$ and a distinguished partial isometry $u$, such that:

\begin{enumerate}
	\item (Cuntz-Krieger relations) $$\sum_{e\in E_1} s_e s_e^* = 1$$ and $$s_g^*s_g=\sum_{e\in E_1:L(i(e))=L(t(g))}s_e s_e^*.$$
	\item For all $e\in E_1 \backslash E_1^{\max}$, $$s_{\lambda (e)}= u s_e \ \text{ and } \  u^*s_{\lambda(e)}=s_e.$$
	\item  $$ u \sum\limits_{e\in E_1^{\max}}s_e = \sum\limits_{e\in E_1^{\min}}s_e u\ \text{ and } \  \sum\limits_{e\in E_1^{\max}}s_e u^* = u^* \sum\limits_{e\in E_1^{\min}}s_e.$$
	\item $u^n$ is a partial isometry, for each $n\in \N$. 
\end{enumerate}

\end{definition}

\begin{obs}
In the case of a Bratteli diagram arising from a primitive, proper substitution, Fujino's C*-algebra $B_{\sigma}$, see \cite{fujino}, is isomorphic to the quotient of $B_E$ by the ideal generated by $1-uu^*$ and $1-u^*u$. 
\end{obs}

Since we have defined the above C*-algebra as an universal object we need to construct a non-zero representation of it. To do so we need to set up notation and give a few definitions.

We define the infinite path space associated to the diagram, denoted by $X$, as the following compact subspace (with the product topology) of $\prod\limits_{n=1}^\infty E_n$: $$X=\left\{\xi\in\prod\limits_{n=1}^\infty E_n:t(\xi_j)=i(\xi_{j+1})\text{ for each  }j\geq 1\right\}$$
and define the finite path space, denoted by $P$, as $$P=\bigcup\limits_{m=1}^\infty\left\{f_1f_2...f_m : f_j\in E_j \text{ and } t(f_j)=i(f_{j+1}) \ \forall j\in \{1,2,...,m-1\}\right\}.$$

Notice that for general Bratteli diagrams there is no clear way to extend the Vershik map. So we will define it on the subsets $X_1:=X\setminus X_{\min}$ and $X_{-1}:=X\setminus X_{\max}$, where
$$X_{\max}=\{\xi\in X:\xi_j\in E_j^{\max},\ \forall j\geq 1\}$$ and $$ \ X_{\min}=\{\xi\in X:\xi_j\in E_j^{\min}, \ \forall j\geq 1\}.$$

Now, to extend the definition of the Vershik map (which so far has only been defined on the edges at a specific level) to  $\lambda:X_{-1}\rightarrow X_1$ we proceed in the following way: given $\xi\in X_{-1}$, let $n$ be such that $\xi_i \in E_i^{\max}$, for $1\leq i \leq n-1$ and $\xi_n\notin E_n^{\max}$ and define $\lambda(\xi)=g_1...g_{n-1}g_n\xi_{n+1}\xi_{n+2}...\in X$, where $g_n$ is the successor of $\xi_n$ (that is, $\lambda(\xi_n)=g_n$) and $g_1...g_{n-1}g_n$ is the unique path such that $g_j\in E_j^{\min}$ for each $j\in\{1,...,n-1\}$. Note that $\lambda:X_{-1}\rightarrow X_1$ is a bijective map. Also notice that the definition of $\lambda$ does not rely on the existence of a unique minimal path (with the lexographic order, as in \cite{DHS}).

In an analogous way we define the successor map, $\lambda$, in a subset of the finite path space. Note that $\lambda:P\setminus P_{\max}\rightarrow P\setminus P_{\min}$, where 
$$P_{\max}=\bigcup\limits_{m=1}^\infty\{f_1...f_m\in P: f_j\in E_j^{\max} \ \forall j =1,..,m \}$$ and $$\ P_{\min}=\bigcup\limits_{m=1}^\infty\{f_1...f_m \in P: f_j\in E_j^{\min}\ \forall j=1,..,m\},$$
is a bijection. 

Finally, for simplicity, from now on we will write $s_g^*s_g= \sum\limits_{e\in E_1:i(e)=t(g)}s_e s_e^*,$ instead of $s_g^*s_g= \sum\limits_{e\in E_1:L(i(e))=L(t(g))}s_e s_e^*$.

We have now set the ground to construct a representation of $B_E$ in the bounded operators in $l^2(X)$, $B(l^2(X))$. To this end we let  $U:l^2(X)\rightarrow l^2(X)$ be defined by $$U(\delta_{\xi})=[\xi\in X_{-1}]\delta_{\lambda(\xi)},$$ where $\{\delta_\xi\}_{\xi\in X}$ is the canonical basis of $l^2(X)$, $[\xi\in X_{-1}]=1$ if $\xi\in X_{-1}$ and $[\xi\in X_{-1}]=0$ if $\xi\notin X_{-1}$. Note that $U$ is a partial isometry and $U^*(\delta_\xi)=[\xi\in X_1]\delta_{\lambda^{-1}(\xi)}$. Also, for each $e\in E_1$, we define $V_e:l^2(X)\rightarrow l^2(X)$ by $$V_e(\delta_\xi)=[t(e)=i(\xi_1)]\delta_{e\xi},$$ which is a partial isometry. Note that the adjoint operator $V_e^*$ of $V_e$ is given by $V_e^*(\delta_\xi)=[\xi_1=e]\delta_{\xi_2\xi_3...}.$

The operators $U$ and $V_e$ satisfy the relations defining $B_E$ and hence, by the universal property of $B_E$, we obtain a *-homomorphism $\phi:B_E\rightarrow B(l^2(X))$ such that $\phi(s_e)=V_e$  and $\phi(u)=U$. Notice that the partial isometry $u\in B_E$ is not an unitary, because $U$ is not an unitary in $l^2(X)$. 

Next we establish that $B_E$ contains a commutative algebra, but to do so we first need to understand better the product structure in $B_E$.


Let $\mathrm W$ be the set of all finite words arising from $(V,E)$, that is, $\mathrm {W} = \displaystyle \bigcup_{m=1}^\infty \{f_1 \ldots f_m: f_j \in E_j \ \forall j=1,..,m \}$. We say that a word $f=f_1 \ldots f_m$ of $\mathrm W$ has length $m$ and write $|f|=m$. Note that a word $f$ is an element of $P$ if and only if $t(f_j)=i(f_{j+1})$ for each $j\in \{1,...,m-1\}$. Finally let $s_f$ denote the corresponding element $s_f=s_{f_1}...s_{f_m}$ in $B_E$.

\begin{lema}\label{lemma1} If $f=f_1...f_m \in \mathrm W$, then $s_f$ is a partial isometry in $B_E$ and:
\begin{enumerate}
\item if $t(f_j)=i(f_{(j+1)})$ for each $j\in\{1,...,m-1\}$, then $s_f\neq 0$;
\item if $t(f_j)\neq i(f_{(j+1)})$ for some $j\in\{1,...,m-1\}$, then $s_f=0$;
\end{enumerate}
\end{lema}

\demo That $s_f$ is a partial isometry follows from the Cuntz-Krieger relations. We prove 1 and 2:

Suppose $t(f_j)=i(f_{j+1})$ for each $j\in \{1,...,m-1\}$ and let $\phi:B_E\rightarrow B(l^2(X))$ be the *-homomorphism defined above. Then
$$\phi(s_f)(\delta_\xi)=V_{f_1}V_{f_2}...V_{f_m}(\delta_\xi)=$$
$$=[t(f_1)=i(f_2)][t(f_2)=i(f_3)]...[t(f_{m-1})=i(f_m)][t(f_m)=i(\xi_1)]\delta_{f_1f_2...f_m\xi}=$$
$$=[t(f_m)=i(\xi_1)]\delta_{f\xi}.$$ It follows that $\phi(s_f)\neq 0$ (since if we choose $\xi\in X$ such that $i(\xi_1)=t(f_m)$ then $\phi(s_f)(\delta_\xi)=\delta_{f\xi}$) and hence $s_f\neq 0$, proving 1.

To prove 2, suppose now that $t(f_j)\neq i(f_{j+1})$ for some $j\in \{1,...,m-1\}$. Then $$s_{f_j}s_{f_{j+1}}=s_{f_j}s_{f_j}^*s_{f_j}s_{f_{j+1}}=s_{f_j}\left(\sum\limits_{e\in E_1:i(e)=t(f_j)}s_es_e^*\right)s_{f_{j+1}}=0,$$ since $s_e^*s_{f_{j+1}}=0$ for each $e\in E_1$ with $i(e)=t(f_j)$ and hence $s_f=s_{f_1}...s_{f_m}=0$.
\fim
  
In the universal C*-algebra $B_E$, for each word $f$, let $e_f=s_fs_f^*$. Then the sub-C*-algebra of $B_E$ generated by $\{e_f:f\in P\}$ is commutative, as we show below:

\begin{proposition}\label{prop1} Let $f,g\in P$, and suppose that $|f|=m$, $|g|=n$, with $m\leq n$. 
\begin{enumerate}
\item If $g=fh$, that is, if $g_1...g_m=f_1...f_m$, then $e_fe_g=e_g=e_ge_f$.
\item If $f_i\neq g_i$ for some $i\in\{1,...,m\}$, then $e_fe_g=0=e_ge_f$.
\end{enumerate}
\end{proposition}

\demo Suppose $g=fh$. Then 
$$e_fe_g=s_fs_f^*s_gs_g^*=s_fs_f^*s_{fh}s_{fh}^*=s_fs_f^*s_fs_hs_{fh}^*=s_fs_hs_{fh}^*=s_{fh}s_{fh}^*=e_g.$$ It follows that $e_fe_g=e_g$, and so $e_ge_f=(e_fe_g)^*=e_g^*=e_g$. This proves 1.

To prove 2, notice first that if $w_1w_2\in P$ then $s_{w_1w_2}^*s_{w_1w_2}=s_{w_1}^*s_{w_1}$, since

$$s_{w_2}^*s_{w_1}^*s_{w_1}s_{w_2}=s_{w_2}^*\left(\sum\limits_{u\in E_1:i(u)=t(w_1)}s_us_u^*\right)s_{w_2}=s_{w_2}^*s_{w_2}.$$

Now, if $f_i\neq g_i$ for some $i\in\{1,...,m\}$, and $k$ is the least index $i$ such that $f_i\neq g_i$, we have that 
$$s_{f_{k}}^*s_{f_{(k-1)}}^*...s_{f_1}^*s_{g_1}...s_{g_{(k-1)}}s_{g_k}=s_{f_{k}}^*s_{f_{(k-1)}}^*...s_{f_1}^*s_{f_1}...s_{f_{(k-1)}}s_{g_k}=s_{f_k}^*s_{f_{(k-1)}}^*s_{f_{(k-1)}}s_{g_k}$$ and $$s_{f_k}^*s_{f_{(k-1)}}^*s_{f_{(k-1)}}s_{g_k}=s_{f_k}^*\left(\sum\limits_{u\in E_1:i(u)=t(f_{(k-1)})}s_us_u^*\right)s_{g_k}=s_{f_k}^*s_{g_k}=0.$$ 

It follows that $s_f^*s_g=0$ and so $e_fe_g=0$ as desired.

\fim

The sub-C*-algebra of $B_E$ generated by $\{e_f:f\in P\}$ will play a crucial role in the next section and so we show below some of its important properties (corresponding, but in a more general setting, to \cite[3.10]{fujino}).


\begin{proposition}\label{prop2} In the universal C*-algebra $B_E$ it holds that:
\begin{enumerate}
\item $\sum\limits_{f\in E_1}e_f=1$;
\item $e_g=\sum\limits_{h\in E_1:i(h)=t(g)}e_{gh}$;
\item $e_{\lambda(f)}=ue_fu^*$ for each $f\in P\setminus P_{\max}$.
\item $e_{\lambda^{-1}(f)}=u^*e_fu$ for each $f\in P\setminus P_{\min}$.
\end{enumerate}
\end{proposition}

\demo 

This proposition is a consequence of the relations which define $B_E$. To give a flavor of the techniques involved we prove 4:

Let $f\in P$, $f=f_1...f_n...f_m$ with $f_i\in E_1^{\min}$ for $i\in \{1,...,n-1\}$ and $f_n\notin E_1^{\min}$. Notice that if $f_k$ and $f_j$ are in $E_1^{\min}$ and $f_k \neq f_j$ then $t(f_k)\neq t(f_j)$. It follows that

$$u^*s_f=u^*s_{f_1}...s_{f_m}=u^*\left(\sum\limits_{e\in E_1^{\min}}s_e\right)s_{f_2}...s_{f_m}=\left(\sum\limits_{e\in E_1^{\max}}s_e\right)u^*s_{f_2}...s_{f_m}=$$
$$=\left(\sum\limits_{e\in E_1^{\max}}s_e\right)...\left(\sum\limits_{e\in E_1^{\max}}s_e\right)u^*s_{f_n}...f_m=$$
$$=\left(\sum\limits_{e\in E_1^{\max}}s_e\right)...\left(\sum\limits_{e\in E_1^{\max}}s_e\right)s_{\lambda^{-1}(f_n)}...s_{f_m}=$$
$$=s_{h_1}...s_{h_{n-1}}s_{\lambda^{-1}(f_n)}...s_{f_m},$$ where $h_1...h_{n-1}\lambda^{-1}(f_n)f_{n+1}...f_m=\lambda^{-1}(f)$. Then $u^*s_f=s_{\lambda^{-1}(f)}$ and hence $e_{\lambda^{-1}(f)}=s_{\lambda^{-1}(f)}s_{\lambda^{-1}(f)}^*=u^*s_fs_f^*u=u^*e_fu.$
\fim

\section{The partial crossed product}

Let $X$ be the infinite path space, $X_{-1}=X\setminus X_{\max}$ and $X_1=X\setminus X_{\min}$, as is the previous section. Recall that $X$ is compact, $X_{-1}, X_1$ are open subsets of $X$ and the Vershik map $\lambda:X_{-1}\rightarrow X_1$ is a homeomorphism. 

Define, inductively, for each $n\in \Z$, $n\geq 2$, $$X_{-n}:=X_{-(n-1)}\cap \lambda^{-(n-1)}(X_{-1}\cap X_{n-1})$$ and $$X_n:=X_{n-1}\cap \lambda^{n-1}(X_1\cap X_{-(n-1)}),$$ which are open subsets of $X$. Let $X_0=X$. Notice that $X_{-n}$ is the domain of $\lambda^n$ and so $\lambda^n$ is a homeomorphism between $X_{-n}$ and $X_n$, for all $n\in \Z$. Now, for each $n \in \Z$, let $\theta_n = \lambda^n$. Then $\theta= \{ \{X_n \}, \{\theta_n\}\}$ is a partial action of $\Z$ in $X$, in the sense of \cite{Ex, MC}, and we denote the corresponding partial action in $C_0(X)$ by $\alpha$. (Note that, for all $n\in \Z$, $\alpha_n:C_0(X_{-n})\rightarrow C_0(X_n)$ is given by $\alpha_n(f)=f\circ \theta_{-n}$). 


Our main goal in the next pages is to show that the partial crossed product $C_0(X)\rtimes_\alpha \Z$ is a sub-C*-algebra of $B_E$. For this purpose, we need to present some notations and preliminary results. 

For each $f=f_1...f_m\in P$, let $1_f$ denote the characteristic function of the subset $X_f=\{\xi\in X:\xi_1...\xi_m=f_1...f_m\}$. Note that $1_f\in C(X)$ because $X_f$ is clopen in $X$. Define $P_0:=P$ and for each $n\in\Z$, $n\geq 1$, $$P_{-n}=\{f\in P:\lambda^{i}(f)\in P\setminus P_{\max}, i=0...n-1\}$$ and $$P_n=\{f\in P:\lambda^{-(i)}(f)\in P\setminus P_{\min}, i=0...n-1 \}.$$ 

Note that if $f=f_1...f_m\in P_n$, for $n\in \Z$, then $1_f\in C_0(X_n)$, since $X_f\subseteq X_n$.

\begin{lema} \label{lemma2} 
\begin{enumerate}
\item Let $f,g\in P$ with $|f|> |g|$. Then $1_f1_g=1_f$, if $f=gh$ and $1_f1_g=0$ otherwise.
\item For each $n\in \Z$, $span \{1_f:f\in P_n\}$ is dense in $C_0(X_n)$.
\item For each $n\in \Z$ and for each $f\in P_{-n}$ it holds that $\alpha_n(1_f)=1_{\lambda^n(f)}$.

\end{enumerate}
\end{lema}

\demo The first statement of the lemma is clear. The second statement follows from the Stone-Weierstrass theorem: $span \{1_f:f\in P_n\}$ vanishes nowhere, since if $\xi \in X_n$ then there exists $m\in \N$ such that $\xi_1...\xi_m\in P_n$, hence $1_{\xi_1...\xi_m}\in C_0(X_n)$ and $1_{\xi_1...\xi_m}(\xi)=1$. Also one can easily check that $span \{1_f:f\in P_n\}$ separates points in $C_0(X_n)$ and the hypothesis of the Stone-Weierstrass theorem are satisfied. Finally, the third statement also holds, since if $f\in P_{-n}$, then $1_f\in C_0(X_{-n})$ and for all $\xi \in X_n$ one has that $$\alpha_n(1_f)(\xi)=1_f(\theta_{-n}(\xi))=1_f(\lambda^{-n}(\xi))=1_{\lambda^n(f)}(\xi).$$
\fim

One important step before we can show that $C(X)\rtimes_\alpha \Z$ is a subalgebra of $B_E$ is to embed $C(X)$ in $B_E$. We do this below.

Let $B$ be the universal C*-algebra generated by projections $\{P_f\}_{f\in P}$ subject to the relations $\sum\limits_{h\in E_1}P_h=1$ and, for each $f\in P$, $P_f=\sum\limits_{h\in E_1:i(h)=t(f)} P_{fh}$. Notice that the first relation implies that all $P_f$, with $f\in E_1$, are orthogonal. This, together with the second relation implies that if $f,g\in P$ with $|f|=|g|$ then $P_f$ and $P_g$ are orthogonal. Finally, this last conclusion, together with the second relation, imply that $B$ is a commutative algebra, with product structure analogous to proposition \ref{prop1}, that is, $P_f P_g = P_g=P_gP_f$ if $g=fh$ and $P_f P_g = 0=P_gP_f$ if $|g|=|f|$ and $g\neq f$ or if $|g|\geq |f|$ and $g \neq fh$. Notice that according to proposition \ref{prop2}, there is a *-representation from $B$ in $B_E$ which maps $P_f$ to $e_f$.

\begin{lema} The spectrum of $B$ is homeomorphic to $X$.
\end{lema}

\demo Let $\widehat{B}$ denote the spectrum of $B$. We will define a map $T:\widehat{B}\rightarrow X$. For this, let $\varphi\in \widehat{B}$. Note that $\varphi(P_f)\in \{0,1\}$, for each $f\in P$, and since $\sum\limits_{f\in E_1}\varphi(P_f)=1$, this implies that there exists one, and only one,  $\xi_1\in E_1$ such that $\varphi(P_{\xi_1})=1$. Now, since $\varphi(P_{\xi_1})=\sum\limits_{h\in E_1}\varphi(P_{\xi_1h})$, there is one, and only one, $\xi_2\in E_1$ such that $\varphi(P_{\xi_1\xi_2})=1$. Continuing in this manner we construct an unique infinite path  $\xi^\varphi\in X$ such that $\varphi(P_{\xi_1^\varphi...\xi_m^\varphi})=1$, for each $m\in \N$. Notice that for $f\in P$, $\varphi(P_f)=1$ if and only if $f=\xi_1^\varphi...\xi_m^\varphi$ for some $m\in \N$. 

So, define $T:\widehat{B}\rightarrow X$ by $T(\varphi)=\xi^\varphi$. 

To see that $T$ is injective, let $\varphi,\psi\in \widehat{B}$ and suppose that $T(\varphi)=T(\psi)$, that is, $\xi^\varphi=\xi^\psi$. Then, for $f\in P$, $\varphi(P_f)=1$ iff $f=\xi_1^\varphi...\xi_m^\varphi=\xi_1^\psi...\xi_m^\psi$ for some $m\in \N$, and this holds iff $\psi(P_f)=1$. It follows that $\varphi$ and $\psi$ coincide on the generators of $B$, and therefore $\varphi=\psi$.

To check that $T$ is surjective, let $\xi\in X$, and define (by the universal property of $B$),  $\varphi:B\rightarrow \C$ by $\varphi(P_f)=1$ if $f=\xi_1...\xi_m$, for some $m\in \N$, and $\varphi(P_f)=0$ otherwise. Then $\varphi\in \widehat{B}$ and $T(\varphi)=\xi$.

It is also not hard to verify that $T$ and $T^{-1}$ are continuous, and hence $T$ is a homeomorphism as desired.
\fim

Let $T:\widehat{B}\rightarrow X$ be the homeomorphism defined in the lema above. It is well known that such a homeomorphism induces a *-isomorphism $\Phi:C(X)\rightarrow C(\widehat{B})$, defined by $\Phi(f)=f\circ T$. Note that for each $f\in P$, $\Phi(1_f)=\widehat{P_f}$, where $\widehat{P_f}$ is the image of $P_f$ by the Gelfand *-isomorphism, since for all $\psi\in \widehat{B}$,  
$$\Phi(1_f)(\psi)=(1_f\circ T)(\psi)=1_f(T(\psi))=1_f(\xi^\psi)=\psi(P_f)=\widehat{P_f}(\psi).$$

So, the composition $\wedge^{-1} \circ \Phi$ is a *-isomorphism from $C(X)$ to $B$ which maps $1_f$ to $P_f$. Since there exists a *-homomorphism from $B$ to $B_E$ that maps $P_f$ to $e_f$, we obtain an *-homomorphism $$\varphi:C(X)\rightarrow B_E$$ such that $\varphi(1_f)=e_f$, for each $f\in P$.

\begin{proposition}\label{lemma4} The *-homomorphism $\varphi:C(X)\rightarrow B_E$ defined above is injective.
\end{proposition}

\demo 
For each $m\geq 1$, let $K_m$ be the cardinality of $\{f\in P:|f|=m\}$ and $C_{K_m}$ be the subalgebra of $C(X)$, generated by $\{1_f:f\in P, |f|=m\}$. Once we show that $\varphi_{|_{C_{K_m}}}:C_{K_m}\rightarrow B_E$ is isometric for all $m\geq 1$, it will follow that $\varphi:C(X)\rightarrow B_E$ is isometric, since any element of $span\{1_f:f\in P\}$ is an element of some subalgebra $C_{K_m}$ (what can be seen from the equality $1_h=\sum\limits_{g\in E_1:i(g)=t(h)}1_{hg}$, for all $h\in P$) and from lema \ref{lemma2} $span\{1_f:f\in P\}$ is dense in $C(X)$.

Now, to show that $\varphi_{|_{C_{K_m}}}$ is isometric, let $B_{K_m}$ be the subalgebra of $B_E$ generated by $\{e_f:f\in P, |f|=m\}$. Let $f^1,f^2,...,f^{K_m}$ be the elements of $\{f\in P:|f|=m\}$. By proposition \ref{prop1}, $B_{K_m}$ is finite dimensional and *-isomorphic to $\C^{d}$, where $d$ is the dimension of $B_{K_m}$. By lema \ref{lemma1} $e_f\neq 0$ for each $f\in P$, and so $d=K_m$. Hence $B_{K_m}$ is *-isomorphic to $\C^{K_m}$, via the *-isomorphism $B_{K_m}\ni e_{f^i}\mapsto e_i \in \C^{K_m}$, where $\{ e_j \}$ is the standard basis of $\C^{K_m}$. Also, $C_{K_m}$ is *-isomorphic to $\C^{K_m}$, via the *-isomorphism $C_{K_m}\ni 1_{f^i}\mapsto e_i \in \C^{K_m}$. So, $C_{K_m}$ is isomorphic to $B_{K_m}$, via the *-isomorphism $C_{K_m}\ni 1_f\mapsto e_f\in B_{K_m}$, and therefore, since $\varphi$ coincides with this isomorphism in $C_{K_m}$, $\varphi_{|_{C_{K_m}}}$ is isometric as desired.
%
\fim

The last step we need before we can embed $C(X)\rtimes_\alpha\Z$ in $B_E$ is to define a partial representation of $\Z$ in $B_E$ that is covariant with respect to the action $\alpha$. For this, define $\pi:\Z\rightarrow B_E$ by $\pi(n)=u^n$ if $n\geq 0$ and $\pi(n)=(u^*)^{|n|}$ if $n<0$. By the definition of $B_E$, $\pi(n)$ is a partial isometry, for each $n\in \Z$. Moreover, we have the following:


\begin{proposition} The map $\pi:\Z\rightarrow B_E$ is a partial representation of $\Z$.
\end{proposition}

\demo
It is clear that $\pi(-n)=\pi(n)^*$ for all $n\in \Z$ and that $\pi(0)=1$. So, all we need to show is that $\pi(-n)\pi(n)\pi(m)=\pi(-n)\pi(n+m)$, for each $n,m\in \Z$. Note that if $m,n$ are both positive or both negative, or if $m=-n$, then $\pi(-n)\pi(n)\pi(m)=\pi(-n)\pi(n+m)$ by the definition of $\pi$. To prove the remaining cases we will use that if $r,s$ are both positive, or both negative, then $\pi(r)\pi(-r)$ and $\pi(-s)\pi(s)$ commute (see \cite[5.3]{Ex1}). So, let $n, m$ be integers with opposite signs.

If $|m|<|n|$ then $n=-m+k$ where $k$ has the same sign as $-m$. So, it follows that  
$$\pi(-n)\pi(n)\pi(m)=\pi(m-k)\pi(k-m)\pi(m)=$$
$$\pi(m)\pi(-k)\pi(k)\pi(-m)\pi(m)=\pi(m)\pi(-m)\pi(m)\pi(-k)\pi(k)=$$ $$=\pi(m)\pi(-k)\pi(k)=\pi(m-k)\pi(k)=\pi(-n)\pi(k)=\pi(-n)\pi(n+m).$$ 

If $|m|>|n|$ then $m=-n+k$, where $k$ has the same sign as $-n$, and in this case, 
$$\pi(-n)\pi(n)\pi(m)=\pi(-n)\pi(n)\pi(-n+k)=\pi(-n)\pi(n)\pi(-n)\pi(k)=$$
$$=\pi(-n)\pi(k)=\pi(-n)\pi(m+n).$$
\fim

\begin{lema}\label{lemma3} For each $n\in \Z$ and $a\in C_0(X_{-n})$ it holds that:
\begin{enumerate}
\item  $\varphi(\alpha_n(a))=\pi(n)\varphi(a)\phi(-n)$ 
\item $\pi(-n)\pi(n)\varphi(a)=\varphi(a)=\varphi(a)\pi(-n)\pi(n)$.
\end{enumerate}
\end{lema}

\demo By lema \ref{lemma2}, it is enough to prove the equalities for $a=1_f$ with $f\in P_{-n}$. By the same lemma, $\alpha_n(1_f)=1_{\lambda^n(f)}$. 

We start with the first statement:

Suppose $n\geq 0$.
Then, $$\pi(n)\varphi(1_f)\pi(-n)=u^ne_f{u^*}^n=e_{\lambda^n(f)}=\varphi(1_{\lambda^n(f)})=\varphi(\alpha_n(1_f)),$$
where the second equality above follows by repeatedly using item 3 of proposition \ref{prop2}. The case $n<0$ follows analogously.

Now, for the second statement, suppose $n\geq 0$. Then, repeatedly using the relations $us_h=s_{\lambda(h)}$ and $u^*s_h=s_{\lambda^{-1}(h)}$ we obtain that $u^ns_f=s_{\lambda^n(f)}$ and ${u^*}^ns_{\lambda^n(f)}=s_f$. It follows that
$$\pi(-n)\pi(n)\varphi(1_f)={u^*}^nu^ne_f={u^*}^nu^ns_fs_f^*=s_fs_f^*=e_f=\varphi(1_f),$$
and taking adjoints we obtain that $\varphi(1_f)\pi(-n)\pi(n)=\varphi(1_f)$.

The case $n<0$ is analogous.
\fim

We are now ready to prove the main result of this section:

\begin{theorem} There exists an injective *-isomorphism from $C(X)\rtimes_\alpha\Z$ in $B_E$.
\end{theorem}

\demo Define $(\varphi\times \pi):\left\{\sum a_n\delta_n:n\in \Z, a_n\in C_0(X_n)\right\}\rightarrow B_E$ by $(\varphi\times \pi)(\sum a_n\delta_n)=\sum \varphi(a_n)\pi(n).$ Note that $\varphi\times \pi$ is contractive, that is, $\|(\varphi\times \pi)(\sum a_n\delta_n)\|\leq \sum \|a_n\|$. Moreover, in light of lema \ref{lemma3} and the fact that $\pi$ is a partial representation, it is straightforward to check that $\varphi\times\pi$ preserves adjoints and multiplications, and hence $\varphi\times \pi$ extends to $C(X)\rtimes_\alpha\Z$. We denote this extension also by $\varphi\times \pi$. Now, by proposition \ref{lemma4}, $\varphi$ is injective and since $\Z$ is amenable and $\theta$ is a free partial action, if follows, by \cite[4.2]{MC} and  \cite[2.6]{Ex3}, that $\varphi\times\pi$ is injective.


\fim

One of the reasons the above result is interesting is that tail equivalence, 
in the case of a simple well ordered Bratteli diagram, played a very important role in the study of Cantor minimal systems and orbit equivalence, see theorem 4.16 of \cite{GPS2} for example. It is easy to check that the the partial orbit equivalence (transformation groupoid) associated to the partial action of the beginning of this section is the same as tail equivalence and hence it is an AF equivalence relation. It follows (by realizing the partial crossed product as the C*-algebra of the partial orbit equivalence relation, see \cite{Abadie} or \cite{Beuter}) that $C(X)\rtimes_\alpha\Z$ is AF. Next we show directly that the partial crossed product $C(X)\rtimes_\alpha\Z$, as in the theorem above, is an AF algebra, regardless of any hypothesis on the Bratteli diagram (this is also a partial converse for a result from Exel, see \cite{Exel}, where he shows that every AF-algebra can be realized as a partial crossed product). Before we proceed we need a technical lemma.

Let $W_N$ be the set of all the words of $P$ with length $N$. Note that $W_N$ is finite, for each $N\in \N$, since $E_1$ is finite.

\begin{lema} For a fixed $N\in \N$, there exists $n_N\in \N$ such that $P_n\cap W_N=\emptyset$, for all $n\in \Z$ with $|n|\geq n_N$.
\end{lema}

\demo Let $m$ be the number of edges in $E_1$. Then $P_{m^N}\cap W_N=\emptyset$ and $P_{-m^N}\cap W_N=\emptyset$. Notice that if $r,s\in \Z$ with $r\geq s\geq 0$ or $r\leq s\leq 0$ then $P_r\subseteq P_s$ and hence it follows that, if $n\in \Z$ with $|n|\geq m^N$ then $P_n\cap W_N=\emptyset$.
\fim

\begin{proposition} The partial crossed product $C(X)\rtimes_\alpha\Z$ associated to the partial action given in the beginning of this section is an AF-algebra.
\end{proposition}

\demo
First notice that $span\{1_f:f\in P_n\}$ is dense in $C_0(X_n)$ (see lemma \ref{lemma2}) and hence $span\{1_f\delta_n:f\in P_n,\,\, n\in \Z\}$ is dense in $C(X)\rtimes_\alpha\Z$. So, it is enough to show that $span \{1_f\delta_n:f\in P_n,\,\,\,n\in \Z\}$ is an increasing union of finite dimensional C*-algebras.

Let $A_N=span\{1_f\delta_n:f\in P_n\cap W_N,\,\,\,n\in \Z\}.$ By the previous lemma, together with the fact that $W_N$ is finite, $A_N$ is a finite dimensional vector space. Next we show that $A_N$ is a finite dimensional C*-algebra.

Let $f\in P_n\cap W_N$ and $g\in P_m\cap W_N$. Then, $$1_f\delta_n1_g\delta_m=\alpha_n(\alpha_{-n}(1_f)1_g)\delta_{n+m}=\alpha_n(1_{\lambda^{-n}(f)}1_g)\delta_{n+m}=[\lambda^{-n}(f)=g]1_f\delta_{n+m}.$$ Notice that if $\lambda^{-n}(f)=g$ then $f\in P_{n+m}$, and so $1_f\delta_{n+m}\in A_N$.
This shows that $A_N$ is closed under multiplication.

To see that $A_N$ is also closed under involution, note that if $f\in P_n$ then $\lambda^{-n}(f)\in P_{-n}$, and hence $(1_f\delta_n)^*=1_{\lambda^{-n}(f)}\delta_{-n}\in A_N$.

So, it follows that $A_N$ is a *-algebra. Since $A_N$ is finite dimensional, $A_N$ is complete and hence a C*-algebra.

It remains to show that $A_N\subseteq A_{N+1}$.

Let $f\in P_n\cap W_N$. 
Notice that we can write $X_f$ as $X_f=\bigcup\limits_{e\in E_1: i(e)=t(f)}^.X_{fe}$, and hence $1_f=\sum\limits_{e\in E_1: i(e)=t(f)}1_{fe}$. Since $f\in P_n\cap W_N$, we have that, for each $e\in E_1$ with $i(e)=t(f)$, $fe\in P_n\cap W_{N+1}$. So $1_f\delta_n=\sum\limits_{e\in E_1:i(e)=t(f)}1_{fe}\delta_n\in A_{N+1}$, and therefore $A_N\subseteq A_{N+1}$.

We conclude that $span\{1_f\delta_n:n\in \Z,\,\,\, f\in P_n\}=\bigcup\limits_{N=1}^\infty A_N$ is an increasing union of finite dimensional C*-algebras and the proposition follows.
\fim

\section{Invariance}

In this last section we show that the C*-algebras introduced in section 2 are almost complete invariants for equivalence of stationary ordered Bratteli diagrams, a notion we make precise below.

Recall that for a Bratteli diagram $(E,V)$, $E_1$ denotes the set of edges at level one, two edges are comparable iff they have the same terminus and $\lambda$ denotes the successor (Vershik) map (either at $E_1$ or at the infinite path space). 

\begin{definition}\label{equivdiagrams} Let $(E,V)$ and $(\widetilde{E},\widetilde{V})$ be two stationary ordered Bratteli diagrams. We say that $(E,V)$ and $(\widetilde{E},\widetilde{V})$ are equivalent if there exists a bijection $T:E_1\rightarrow \widetilde{E}_1$ such that:
\begin{enumerate}
\item For all $e,f \in E_1$, $t(e)=i(f)$ if and only if $t(T(e))=i(T(f))$;
\item For all $e \in E_1\setminus E_1^{\max}$, $T(\lambda(e))=\widetilde{\lambda}(T(e))$, where we assume that $\lambda$ and $\widetilde{\lambda}$ are defined only on $E_1\setminus E_1^{\max}$ and $\widetilde{E_1}\setminus \widetilde{E}_1^{\max}$, respectively.
\end{enumerate}
\end{definition}

\begin{obs} It follows from the second item in the above definition that $e$ is maximal if and only if $T(e)$ is maximal, and $e$ is minimal if and only if $T(e)$ is minimal.
\end{obs}

\begin{example} Let $(E,V)$ end $(\widetilde{E},\widetilde{V})$ be two Bratteli diagrams as in the picture below. In this case, $E_1=\{e_1, e_2, e_3, e_4, e_5, e_6\}$ with $e_1<e_2$ and $e_4<e_5<e_6$, and $\widetilde{E}_1=\{f_1, f_2, f_3, f_4, f_5, f_6\}$ with $f_1<f_2$ and $f_4<f_5<f_6$. 
The Bratteli diagrams $E$ and $\widetilde{E}$ are equivalent, with the bijection $T:E\rightarrow \widetilde{E}$ given by by $T(e_n)=f_n$, for $n\in\{1,2,...,6\}$.  
\end{example}

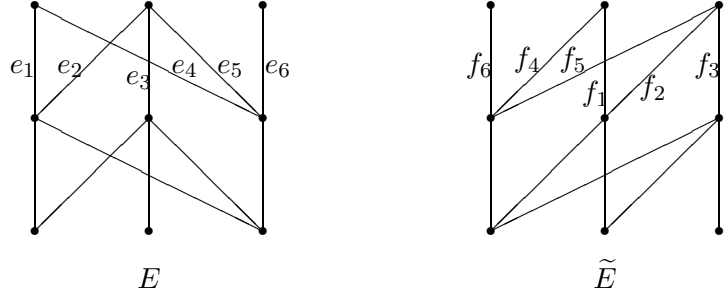
\centerline{
\setlength{\unitlength}{1.5cm}
\begin{picture}(6,0.6)
\put(0,0){\circle*{0.08}}
\put(1,0){\circle*{0.08}}
\put(2,0){\circle*{0.08}}
\put(0,-1){\circle*{0.08}}
\put(1,-1){\circle*{0.08}}
\put(2,-1){\circle*{0.08}}
\put(0,-2){\circle*{0.08}}
\put(1,-2){\circle*{0.08}}
\put(2,-2){\circle*{0.08}}
\put(0,0){\line(2,-1){2}}
\put(1,0){\line(1,-1){1}}
\put(2,0){\line(0,-1){1}}
\put(0,-1){\line(2,-1){2}}
\put(1,-1){\line(1,-1){1}}
\put(2,-1){\line(0,-1){1}}
\put(0,0){\line(0,-1){1}}
\put(1,0){\line(0,-1){1}}
\put(1,0){\line(-1,-1){1}}
\put(0,-1){\line(0,-1){1}}
\put(1,-1){\line(0,-1){1}}
\put(1,-1){\line(-1,-1){1}}
\put(-0.22,-0.6){$e_1$}
\put(0.2,-0.6){$e_2$}
\put(0.8,-0.7){$e_3$}
\put(1.2,-0.6){$e_4$}
\put(1.6,-0.6){$e_5$}
\put(2.02,-0.6){$e_6$}
\put(4,0){\circle*{0.08}}
\put(5,0){\circle*{0.08}}
\put(6,0){\circle*{0.08}}
\put(4,-1){\circle*{0.08}}
\put(5,-1){\circle*{0.08}}
\put(6,-1){\circle*{0.08}}
\put(4,-2){\circle*{0.08}}
\put(5,-2){\circle*{0.08}}
\put(6,-2){\circle*{0.08}}
\put(4,0){\line(0,-1){1}}
\put(5,0){\line(0,-1){1}}
\put(6,0){\line(0,-1){1}}
\put(4,-1){\line(0,-1){1}}
\put(5,-1){\line(0,-1){1}}
\put(6,-1){\line(0,-1){1}}
\put(5,0){\line(-1,-1){1}}
\put(5,-1){\line(-1,-1){1}}
\put(6,0){\line(-1,-1){1}}
\put(6,-1){\line(-1,-1){1}}
\put(6,0){\line(-2,-1){2}}
\put(6,-1){\line(-2,-1){2}}
\put(3.78,-0.6){$f_6$}
\put(4.2,-0.55){$f_4$}
\put(4.6,-0.55){$f_5$}
\put(4.8,-0.85){$f_1$}
\put(5.3,-0.8){$f_2$}
\put(5.78,-0.6){$f_3$}
\put(0.9,-2.5){$E$}
\put(4.9,-2.5){$\widetilde{E}$}
\end{picture}}

\vspace{3.6cm}

\begin{theorem}\label{teorIsom} If $(E,V)$ and $(\widetilde{E},\widetilde{V})$ are equivalent then $B_E$ and $B_{\widetilde{E}}$ are isomorphic C*-algebras. Furthermore, there exists an *-isomorphism $\phi:B_E\rightarrow B_{\widetilde{E}}$ that preserves the generators of the algebras, that is, $\phi$ is such that $\{\phi(s_e)\}_{e\in E_1} =\{s_{\widetilde{e}}\}_{\widetilde{e}\in \widetilde{E_1}}$ and $\phi(u)=\widetilde{u}$, where $\{s_e\}_{e\in E_1}\cup \{u\}$ are the generators of $B_E$ and $\{s_{\widetilde{e}}\}_{\widetilde{e}\in \widetilde{E}_1}\cup \{\widetilde{u}\}$ are the generators of $B_{\widetilde{E}}$.
\end{theorem}

\demo
Let $\{s_e\}_{e\in E_1}$ and $u$ be the generators of $B_E$ and $\{s_{\widetilde{e}}\}_{\widetilde{e}\in \widetilde{E}_1}$ and $\widetilde{u}$ be the generators of $B_{\widetilde{E}}$. By the universal property of $B_E$ there exists a *-homomorphism $\phi:B_E\rightarrow B_{\widetilde{E}}$ such that $\phi(s_e)=s_{T(e)}$ and $\phi(u)=\widetilde{u}$. In a similar way there is a *-homomorphism $\psi:B_{\widetilde{E}}\rightarrow B_E$ such that $\psi(\widetilde{u})=u$ and $\psi(s_{\widetilde{e}})=s_{T^{-1}(\widetilde{e})}$. It is now straightforward to check that $\psi$ is the inverse of $\phi$.
\fim

\begin{theorem}
Let $(E,V)$ and $(\widetilde{E},\widetilde{V})$ be two stationary ordered Bratteli diagrams. If there exists a *-isomorphism $\psi:B_E\rightarrow B_{\widetilde{E}}$ that preserves the generators of the algebras (as in theorem \ref{teorIsom}), then $(E,V)$ and $(\widetilde{E}, \widetilde{V})$ are equivalent.
\end{theorem}

\demo 

Define $T:E\rightarrow \widetilde{E}$ by $T(e)=\widetilde{e}$, where $\widetilde{e}$ is such that $\psi(s_e)=s_{\widetilde{e}}$. We have to check that $T$ is a bijection that satisfies conditions 1 and 2 of definition \ref{equivdiagrams}. For what follows it is important to notice that $psi(s_e)=s_{T(e)}$ for all $e\in E_1$. 

That $T$ is surjective is clear. Also, if $e\neq f$ in $E_1$, then $s_{T(e)}=\psi(s_e)\neq \psi(s_f)=s_{T(f)}$, and so $T(e)\neq T(f)$. Therefore $T$ is a bijection. 

Next we verify the first item of definition \ref{equivdiagrams}. Let $e\in E_1$. Then $$\sum\limits_{i(g)=t(T(e))}s_gs_g^*=s_{T(e)}^*s_{T(e)}=\psi(s_e^*s_e)=\psi\left(\sum\limits_{i(f)=t(e)}s_fs_f^*\right)=\sum\limits_{i(f)=t(e)}s_{T(f)}s_{T(f)}^*$$
and so $$\sum\limits_{i(g)=t(T(e))}s_gs_g^*=\sum\limits_{i(f)=t(e)}s_{T(f)}s_{T(f)}^*.$$

Now, for a fixed $g_0\in \{g\in \widetilde{E_1}: i(g)=t(T(e))\}$ we have that $$s_{g_0}=\sum\limits_{i(g)=t(T(e))}s_gs_g^*s_{g_0}=\sum\limits_{i(f)=t(e)}s_{T(f)}s_{T(f)}^*s_{g_0},$$ and hence there exists $T(f)\in \{T(f):i(f)=t(e)\}$ such that $T(f)=g_0$. 

We have shown that $$\{g\in \widetilde{E_1}:i(g)=t(T(e))\}\subseteq \{T(f):i(f)=t(e)\},$$ and proceeding is a similar way we also obtain that $$\{T(f):i(f)=t(e)\}\subseteq\{g\in \widetilde{E_1}:i(g)=t(T(e))\},$$ and hence $$\{T(f):i(f)=t(e)\}=\{g\in \widetilde{E_1}:i(g)=t(T(e))\}.$$

So, to prove condition 1 of definition \ref{equivdiagrams}, let $e,f \in E_1$ be such that $i(f)=t(e)$. Then $T(f)=g$ for some $g\in \widetilde{E_1}$ such that $i(g)=t(T(e))$, and hence $i(T(f))=i(g)=t(T(e))$. On the other hand, if $i(T(h))=t(T(e))$, for $h,e \in E_1$,  then $T(h)\in \{g\in \widetilde{E_1}:i(g)=t(T(e))\}=\{T(f):i(f)=t(e)\}$, and so $i(h)=t(e)$.

Before we verify the second item in definition \ref{equivdiagrams}, we need to show that if $e\in E_1$ is not a maximal edge then $T(e)$ is not a maximal edge either.

Let $e\in E_1$ be a non maximal edge and suppose $T(e)$ is maximal. Note that $\widetilde{u}s_{T(e)}=\psi(us_e)=\psi(s_{\lambda(e)})=s_{T(\lambda(e))}$. In terms of the representation of the algebra introduced in section 2, this equality yields the equality $UV_{T(e)}=V_{T(\lambda(e))}$, that is, $UV_{T(e)}(\delta_{\xi})=V_{T(\lambda(e))}(\delta_\xi)$ for each $\xi$ in the path space $X$ associated to $\widetilde{E}$.

Notice that there exists $\xi \in X$, non maximal, such that $T(e)\xi\in X$, since if we suppose that no such $\xi$ exists then $UV_{T(e)}=0$, what implies that $V_{T(\lambda(e))}=0$, a contradiction (since $V_f\neq 0$ for each $f\in \widetilde{E}_1$). 
So, for this $\xi$, we have that $\delta_{\widetilde{\lambda}(T(e)\xi)}=UV_{T(e)}\delta_\xi=V_{T(\lambda(e))}(\delta_\xi)=\delta_{T(\lambda(e))\xi}$, and it follows that $\widetilde{\lambda}(T(e)\xi)=T(\lambda(e))\xi$. Since $T(e)$ is maximal then $\widetilde{\lambda}(T(e)\xi)=f\widetilde{\lambda}(\xi)$, where $f\in \widetilde{E_1}$, and therefore $\widetilde{\lambda}(\xi)=\xi$, a contradiction. 
It follows that $T(e)$ is not maximal. 

Let us finally verify the second item of definition \ref{equivdiagrams}. Let $e$ be a non maximal edge in $E_1$. Then $T(e)$ is not maximal and $$s_{\widetilde{\lambda}(T(e))}=\widetilde{u}s_{T(e)}=\psi(us_e)=\psi(s_{\lambda(e)})=s_{T(\lambda(e))},$$ and so $T(\lambda(e))=\widetilde{\lambda}(T(e))$ as desired.
\fim

\vspace{1.5pc}

Daniel Gonçalves, Departamento de Matemática, Universidade Federal de Santa Catarina, Florianópolis, 88040-900, Brasil

Email: daemig@gmail.com

\vspace{0.5pc}
Danilo Royer, Departamento de Matemática, Universidade Federal de Santa Catarina, Florianópolis, 88040-900, Brasil

Email: royer@mtm.ufsc.br
\vspace{0.5pc}


\addcontentsline{toc}{section}{References}
\begin{thebibliography}{99}

\bibitem{Abadie} F. Abadie, {\it On partial actions and groupoids}, Proc. Am. Math. Soc., 132(4) (2004), 1037- 1047.

\bibitem{Beuter} V. Beuter, D. Gonçalves, {\it Algebras from Equivalence Relations}, Submitted for publication.

\bibitem{DHS} F. Durand, B. Host and C. Skau, {\it Substitutional dynamical systems, Bratteli diagrams and dimension groups.} Ergodic Theory and Dynamical Systems, 19 (1999) , pp 953-993.

\bibitem{Exel} R. Exel, {\it Approximately finite C*-algebras and partial automorphisms}, Math. Scand., 77 (1995), 281 - 288.


\bibitem{Ex} R. Exel, {\it Circle actions on C*-algebras, partial automorphisms and a generalized Pimsner-Voiculescu exact sequence}, J. Funct. Analysis, 122 (1994), 361 - 401.

\bibitem{Ex1} R. Exel, {\it Partial representations and amenable Fell bundles over free groups}, Pacific J. Math., 192 (2000), 39-63.

\bibitem{Ex3} R. Exel, M. Laca and J. Quigg {\it Partial dynamical systems and C*-algebras generated by partial isometries},	J. Operator Theory, 47 (2002), 169 - 186.

\bibitem{fujino} M. Fujino, {\it C*-algebras arising from substitutions}, Ergodic Theory and Dynamical Systems, 30 (2010), 1685-1702. 


\bibitem{HPS} R.H. Herman, I. Putnam, C. Skau, {\it Ordered Bratteli diagram, dimension groups, and topological dynamics}, Int. J. Math. 6 (1992), 827-864.

\bibitem{MC} K. McClanahan, {\it K-theory for partial actions by discrete groups}, J. Funct. Anal., 130 (1995), 77-117.


\bibitem{GPS2} T. Giordano, C.F. Skau and I. Putnam, {\it Affable equivalence relations and orbit structure of Cantor dynamical systems}, Ergodic Theory and Dynamical Systems, 24 (2004), 441-476.

\bibitem{PY} Pivato, Yassawi, {\it Embedding Bratteli-Vershik systems in cellular
automata}, Ergodic Theory and Dynamical Systems, 30 (2010), 1561-1572.

\bibitem{Putnam} I. Putnam, {\it The C*-algebras associated with minimal homeomorphisms of the Cantor set}, Pacific J. Math., 136(2) (1989), 329-353. 



\end{thebibliography}
\end{document}